\documentclass{elsart}

\usepackage{epsfig}

\begin{document}

\begin{frontmatter}

\title{Zeros of the Macdonald function of complex order}

\author{Erasmo M. Ferreira}
\address{Instituto de F\'{\i}sica, Universidade Federal do Rio de
Janeiro, Caixa Postal 68528, 21941-972 Rio de Janeiro, RJ, Brasil}

\author{Javier Sesma}
\address{Departamento de F\'{\i}sica Te\'orica,
         Facultad de Ciencias, Universidad de Zaragoza, 50009 Zaragoza, Spain}

\begin{abstract}
The $z$-zeros of the modified Bessel function of the third kind
$K_{\nu}(z)$, also known as modified Hankel function or Macdonald
function, are considered for arbitrary complex values of the order
$\nu$. Approximate expressions for the zeros, applicable in the
cases of very small or very large $|\nu|$, are given. The
behaviour of the zeros for varying $|\nu|$ or $\arg \nu$, obtained
numerically, is illustrated by means of some graphics.
\end{abstract}


\end{frontmatter}

\begin{tabbing}
{\it Corresponding author:} \quad \=  Javier Sesma \\
{\it address:} \> Departamento de Fisica Teorica, \+ \\
                  Facultad de Ciencias,\\
                  50009 Zaragoza, Spain. \- \\
{\it phone:} \> 34 - 976 761 265 \\
{\it fax:}   \> 34 - 976 761 264 \\
{\it e-mail:} \> javier@unizar.es
\end{tabbing}

\bigskip

\section{Introduction.}

The relevance of Bessel functions to the solution of a great
variety of problems in Physics and Engineering is widely known and
does not need to be stressed. In particular, zeros of the
different kinds of those functions are closely related to energies
of bound or resonant physical systems. For instance, the
$\nu$-zeros of the Hankel function $H_{\nu}(z)$ for positive
values of $z$ determine the poles of the amplitude of scattering
of various kinds of waves by spheres and cylinders \cite{kell}.
Almost four decades ago we considered \cite{ferr} the zeros of the
modified Bessel function of the third kind, $K_{\nu}(z)$, also
known as modified Hankel function or Macdonald function.
Information about the $z$-zeros of $K_{\nu}(z)$ for positive order
$\nu$ was found in the book by Watson \cite{wats}. However, our
interest, motivated by a quantum mechanical problem, was on
$z$-zeros when $\nu$ assumes purely imaginary values. For that
case it was already known \cite{frie,gray} that $K_{\nu}(z)$
presents an infinity of zeros, all of them located on the positive
real semiaxis of the complex $z$-plane. We showed that these
positive zeros, labelled according to decreasing values, form a
sequence tending to be a geometric progression with ratio
$\exp(-\pi/|\nu|)$. We gave, also, asymptotic approximations to
the position of the largest zeros for large values of $|\nu|$.
Further discussions of the properties of the zeros of
$K_{\nu}(z)$, $\nu$ pure imaginary, have been carried out by
Laforgia \cite{lafo} and by Dunster \cite{duns}. Zeros of
$K_{\nu}(z)$ with complex $\nu$ have been less studied. In view of
their application in the analysis of electromagnetic wave
propagation in anisotropic inhomogeneous conducting media, Nalesso
\cite{nale} has considered the $\nu$-zeros of $K_{\nu}(t\nu)$ for
fixed positive values of $t$. But, to our knowledge, a discussion
of the $z$-zeros of $K_{\nu}(z)$ for arbitrary complex $\nu$ is
still lacking.

In recent years, stimulated by review papers by Lozier and Olver
\cite{loz1,loz2} and by the DLMF project, a considerable number of
new algorithms for the computation of special functions have been
published. Concerning the Macdonald function, Refs.
\cite{fab2,gaut,gil1,gil2,gil3,gil4,segu} are to be noticed. As a
general rule, algorithms lose accuracy in the vicinity of the
zeros. Hence the interest of having a knowledge, at least
approximate, of the location of those zeros. On the other hand,
$K_{\nu}(z)$ presents only a finite number of $z$-zeros in the
Riemann sheet $-\pi<\arg z\leq\pi$ if $\nu$ is real, whereas an
infinity of zeros occur if $\nu$ is pure imaginary. It has
therefore seemed to us interesting to study the evolution of the
$z$-zeros of $K_{\nu}(z)$ as the modulus and/or argument of $\nu$
are continuously changed.

In what follows, we consider only values of $\nu$ in the quadrant
\[ 0\leq\arg \nu\leq\pi/2,\]
and $z$ is assumed to lie in the principal Riemann sheet,
$-\pi<\arg z\leq\pi$, for large values of $\nu$. The symmetry
relations
\begin{equation}
K_{-\nu}(z)=K_{\nu}(z),\qquad
K_{\overline{\nu}}(\overline{z})=\overline{K_{\nu}(z)},
\label{uno}
\end{equation}
allow to extend our results to values of $\nu$ in all quadrants.
On the other hand, the relation \cite[Eq. 9.6.4]{abra}
\begin{equation}
K_{\nu}(z)={\textstyle{\frac{1}{2}}}\pi i e^{\frac{1}{2}\nu\pi
i}H_{\nu}^{(1)}\left(ze^{\frac{1}{2}\pi i}\right), \label{dos}
\end{equation}
valid when $K$ and $H^{(1)}$ are analytically continued to every
Riemann sheet, proves that the pattern of $z$-zeros of the Hankel
function $H_{\nu}^{(1)}(z)$ results from that of $K_{\nu}(z)$ by
rotation by an angle of $+\pi/2$ around the origin.

We discuss in Section 2 the behaviour of the zeros for decreasing
values of $|\nu|$. The case of large values of $|\nu|$ is
considered in Section 3, where asymptotic approximations to the
zeros are given. Numerical results for the zeros of $K_{\nu}(z)$
with moderate $|\nu|$ are presented graphically in Section 4.
Finally, the particular case of real $\nu$ is commented in Section
5.

\section{Zeros of $K_{\nu}(z)$ for $|\nu|\ll1$.}

The particular case of pure imaginary $\nu$ revealed \cite{ferr}
that, as $|\nu|$ goes to zero, all zeros of $K_{\nu}(z)$ approach
the origin. Therefore we start by assuming $|z|\ll1$ in the
expression of $K_{\nu}(z)$ as a sum of two convergent ascending
series \cite[Eqs. 9.6.2 and 9.6.10]{abra}
\begin{equation}
K_{\nu}(z) = \frac{(\pi/2)}{\sin(\nu\pi)}\left(
\frac{(z/2)^{-\nu}}{\Gamma(1\!-\!\nu)}\sum_{k=0}^{\infty}
\frac{(z^2/4)^k}{k!(1\!-\!\nu)_k} -
\frac{(z/2)^{\nu}}{\Gamma(1\!+\!\nu)}\sum_{k=0}^{\infty}
\frac{(z^2/4)^k}{k!(1\!+\!\nu)_k} \right),  \label{tres}
\end{equation}
valid for $\nu$ different from an integer. Obviously, the zeros of
$K_{\nu}(z)$ should satisfy the relation
\begin{equation}
(z/2)^{2\nu} = \frac{\Gamma(1\!+\!\nu)}{\Gamma(1\!-\!\nu)}\;
\frac{{\displaystyle{\sum_{k=0}^{\infty}
\frac{(z^2/4)^k}{k!(1\!-\!\nu)_k}}}}{{\displaystyle{\sum_{k=0}^{\infty}
\frac{(z^2/4)^k}{k!(1\!+\!\nu)_k}}}}, \label{cuatro}
\end{equation}
that, in the case of being $|z|\ll1$, can be approximated by
\begin{equation}
(z/2)^{2\nu} = \frac{\Gamma(1\!+\!\nu)}{\Gamma(1\!-\!\nu)}
\left(1+\frac{\nu}{2(1\!-\!\nu^2)}\,z^2+O(z^4)\right), \qquad
|z|\ll1. \label{cinco}
\end{equation}
By taking logarithms, one obtains, as a first approximation,
\begin{equation}
\hspace{-.5cm}\log |z| + i\arg z \simeq -\frac{n\pi i}{\nu} + \log
2 +\frac{1}{2\nu}\log (\Gamma(1\!+\!\nu)/\Gamma(1\!-\!\nu))
+\frac{z^2}{4(1\!-\!\nu^2)},  \label{seis}
\end{equation}
where $n$ could in principle take the values $n=0, \pm 1, \pm 2,
\ldots$, although, according to that obtained in the case of $\nu$
being pure imaginary, actually existing zeros correspond to
\[ n = 1, 2, 3, \ldots .  \]
These values of $n$ will be used as a label to characterize each
zero in the form $z_n$. Ignoring, in the crudest approximation,
the last term in the right hand side of Eq. (\ref{seis}), one has
\begin{eqnarray}
\log |z_n| & \simeq & \Re \left(\frac{1}{\nu}(-n\pi i +
{\textstyle{\frac{1}{2}}}\log
(\Gamma(1\!+\!\nu)/\Gamma(1\!-\!\nu))\right)+\log 2,   \label{siete} \\
\arg z_n & \simeq & \Im \left(\frac{1}{\nu}(-n\pi i +
{\textstyle{\frac{1}{2}}}\log
(\Gamma(1\!+\!\nu)/\Gamma(1\!-\!\nu))\right). \label{ocho}
\end{eqnarray}
The behavior of the zeros $z_n$ as $|\nu|$ decreases can be
immediately obtained from these expressions by using the series
expansion \cite[Eq. 6.1.34]{abra}
\begin{equation}
\frac{1}{\Gamma(1\!+\!\nu)} = 1+\gamma\nu+c_3\nu^2+\ldots
\label{nueve}
\end{equation}
to obtain
\begin{equation}
\frac{\Gamma(1\!+\!\nu)}{\Gamma(1\!-\!\nu)} = 1 - 2\gamma\nu
+2\gamma^2\nu^2+\ldots, \label{diez}
\end{equation}
where $\gamma$ represents the Euler constant, $\gamma =
0.5772156649\ldots$. By retaining terms up to $O(|\nu|^{-1})$, it
turns out
\begin{eqnarray}
\log |z_n| & \simeq & -n\pi\,\frac{\Im \nu}{|\nu|^2} + \log 2 -
\gamma, \label{duno}
\\ \arg z_n & \simeq & -n\pi\,\frac{\Re \nu}{|\nu|^2}. \label{ddos}
\end{eqnarray}
These two equations allow us to obtain a picture of the behaviour
of the zeros of $K_{\nu}(z)$ as $|\nu|$ tends to zero. For
instance, let us assume that $|\nu|$ decreases with $\arg \nu$
fixed. If $\arg \nu = \pi/2$, the zeros $z_n$ approach the origin
along the positive real semiaxis, as described in Ref.
\cite{ferr}. If $0 < \arg \nu < \pi/2$, $|z_n|$ tends to zero
whereas $\arg z_n$ decreases without limit. In other words, the
zeros approach the origin spiraling clockwise infinitely.
According to Eq. (\ref{ddos}), most of zeros (or even all of them,
for sufficiently small $|\nu|$) occur outside the principal
Riemann sheet. For $\arg \nu = 0$, Eqs. (\ref{duno}) and
(\ref{ddos}) are no more valid because, as we will see later, the
possibility exists of $z_n$ going to infinity for certain values
of $\nu$. The behaviour of the zeros in this case is much more
complicated and will be described in Section 5.

\section{Zeros of $K_{\nu}(z)$ for $|\nu|\gg 1$.}

We have already mentioned that Eq. (\ref{dos}) allows one to
obtain the zeros of the Macdonald function from those of the
Hankel function, and viceversa. Previous research, by other
authors \cite{coch,fran,heth,kell,magn,sand}, on the $\nu$-zeros
of the Hankel function shows that, for large values of $|z|$, they
occur at $\nu \simeq z$. Asymptotic expansions of $H^{(1)}$ valid
in the transition region should then be used to obtain approximate
expressions of its $z$-zeros  for $|\nu|\gg 1$. From Eqs. 9.3.23
and 9.3.24 of Ref. \cite{abra}, by using Eqs. 9.1.3 and 10.4.9, a
convenient expansion can be written, namely
\begin{eqnarray}
\lefteqn{H_{\nu}^{(1)}(\nu+\theta\nu^{1/3}) \sim
\frac{2^{4/3}}{\nu^{1/3}}\,e^{-i\pi/3} \Bigg(
{\mbox{Ai}}(-2^{1/3}e^{i2\pi/3}\theta)
\sum_{k=0}^{\infty}f_k(\theta)\nu^{-2k/3}} \nonumber \\
& & \hspace{3.5cm} +\frac{2^{1/3}}{\nu^{2/3}}\,e^{i2\pi/3}
{\mbox{Ai}}^{\prime}(-2^{1/3}e^{i2\pi/3}\theta)
\sum_{k=0}^{\infty}g_k(\theta)\nu^{-2k/3}\Bigg), \label{dtres}
\end{eqnarray}
where $f_0(\theta)=1$, the other functions $f_k(\theta)$ and
$g_k(\theta)$ being given in Eqs. 9.3.25 and 9.3.26 of the same
reference. A more powerful asymptotic expansion of $H^{(1)}$
exists \cite[Eq. 9.3.37]{abra}, but its coefficients are much more
complicated and, although reversion of that expansion could be
done by means of a procedure due to Fabijonas and Olver
\cite{fab1}, we find preferable, for the sake of simplicity and
transparency, to work with Eq. (\ref{dtres}). A simple inspection
of that expansion shows that, for large $|\nu|$, it vanishes for
values of $\theta$ verifying
\begin{equation}
\theta_s = -2^{-1/3}\,e^{-i2\pi/3}\,(a_s+\varepsilon_s(\nu)),
\qquad s=1, 2, 3, \ldots, \label{dcuatro}
\end{equation}
where $a_s$ represents each one of the zeros of the Airy function,
given in Table 10.13 of Ref. \cite{abra}, and
$\varepsilon_s(\nu)=O(\nu^{-2/3})$ is to be determined. Trying an
expansion of the form
\begin{equation}
\varepsilon_s(\nu) = \sum_{j=1}^{\infty} \alpha_{s,j} \,
\nu^{-2j/3}, \label{dcinco}
\end{equation}
and using Taylor expansions of the Airy function and of its
derivative around $a_s$, it is straightforward, with the aid of a
computer algebra software package, to check that the right hand
side of Eq. (\ref{dtres}) cancels out for
\begin{eqnarray}
\alpha_{s,1} & = & -\frac{3}{10}\,2^{-1/3}\,e^{-i2\pi/3}\,a_s^2, \nonumber \\
\alpha_{s,2} & = & -\frac{1}{700}\,2^{1/3}\,e^{i2\pi/3}\,(a_s^3+10), \nonumber\\
\alpha_{s,3} & = & \frac{1}{126000}\,a_s\,(479\,a_s^3-40), \nonumber\\
\alpha_{s,4} & = & \frac{1}{16170000}\,2^{-1/3}\,e^{-i2\pi/3}\,
a_s^2\,(20231\,a_s^3+55100), \qquad \mbox{etc}.  \nonumber
\end{eqnarray}
Approximate values of the $z$-zeros of $H_{\nu}^{(1)}(z)$ are then
given by
\begin{equation}
z_s \sim (\nu+\theta_s\,\nu^{1/3}), \qquad |\nu|\gg 1, \qquad s=1,
2, 3, \ldots. \label{dseis}
\end{equation}
Reversion of this relation between $z$ and $\nu$ would provide an
approximation for the $\nu$-zeros of $H_{\nu}^{(1)}(z)$ when $|z|$
is large,
\begin{eqnarray}
\lefteqn{\hspace{-.5cm}\nu_s \sim z+2^{-1/3}e^{-i2\pi/3}a_s\,
z^{1/3}+\frac{1}{60}2^{1/3}e^{i2\pi/3}
a_s^2\,z^{-1/3}-\frac{1}{700}(a_s^3+10)\,z^{-1}} \nonumber \\
 & & \hspace{.6cm}+ \frac{1}{1134000}2^{-1/3}e^{-i2\pi/3}a_s\,
(281\,a_s^3+10440)\,z^{-5/3} \nonumber \\
 & & \hspace{.6cm} - \frac{1}{2619540000}2^{1/3}e^{i2\pi/3}a_s^2\,
(73769\,a_s^3+6624900)\,z^{-7/3}, \quad |z|\gg 1, \label{dsiete}
\end{eqnarray}
to be compared with those given in Refs.
\cite{coch,fran,heth,kell,sand}. According to Eq. (\ref{dos}), the
$z$-zeros of $K_{\nu}^{(1)}(z)$ are given approximately by
\begin{equation}
z_s \sim e^{-i\pi/2}\,(\nu+\theta_s\,\nu^{1/3}), \qquad |\nu|\gg
1, \qquad s=1, 2, 3, \ldots. \label{docho}
\end{equation}

The usefulness of Eqs. (\ref{dsiete}) and (\ref{docho}) is
restricted to the first values of $s=1,2,3,\ldots$, due to the
fact that $|a_s|$ increases with $s$. This makes the the omitted
terms in the expansion (\ref{dsiete}) to increase and invalidates
their suppression. Besides this, $|\theta_s|$, given by
(\ref{dcuatro}), increases with $s$ and, for a given $\nu$, may
reach a value such that expansion (\ref{dtres}) ceases from being
useful. If zeros corresponding to a large $s$ are to be
considered, approximate values can be obtained following a
procedure, already used by Magnus and Kotin \cite{magn} and by
Cochran \cite{coch}, consisting in taking logarithms in Eq.
(\ref{cuatro}),
\begin{eqnarray}
\lefteqn{\hspace{-.5cm} 2\,\nu\log(z/2)+i\,2n\pi
=\log(\Gamma(\nu))-\log(-\Gamma(-\nu))} \nonumber \\
 & & +\log\left(\sum_{k=0}^{\infty}
\frac{(z^2/4)^k}{k!(1\!-\!\nu)_k}\right)-\log\left(\sum_{k=0}^{\infty}
\frac{(z^2/4)^k}{k!(1\!+\!\nu)_k}\right), \quad n=0,\pm 1,\pm 2,
\ldots , \label{dnueve}
\end{eqnarray}
and approximating the logarithms of the gamma functions by their
asymptotic expansions \cite[Eq. 6.1.40]{abra}. One obtains in this
way
\begin{eqnarray}
\lefteqn{\log(z/2)\sim\log\nu -1 -i\pi/2
-i(n-{\textstyle{\frac{1}{4}}})\pi/\nu  +
\sum_{m=1}^{\infty}\frac{B_{2m}}{2m(2m-1)\nu^{2m}}}\nonumber \\
 & & \hspace{2cm} + \frac{1}{2\nu}\left(\log\left(\sum_{k=0}^{\infty}
\frac{(z^2/4)^k}{k!(1\!-\!\nu)_k}\right)-\log\left(\sum_{k=0}^{\infty}
\frac{(z^2/4)^k}{k!(1\!+\!\nu)_k}\right)\right). \label{veinte}
\end{eqnarray}
Assuming $|z|\ll|\nu|$ and keeping the dominant terms, it comes
out for the $z$-zeros of $K_{\nu}(z)$
\begin{equation}
z_n \simeq
e^{-i\pi/2}\,2\,\nu\,\exp\left(-1-i(n-{\textstyle{\frac{1}{4}}})
\pi/\nu\right). \label{vuno}
\end{equation}
Obviously, since we are considering values of $\nu$ in the first
quadrant, only large positive values of the integer $n$ are
compatible with the assumption $|z|\ll|\nu|$. Moreover, for a
given value of $|\nu|$, the approximation (\ref{vuno}) looses
accuracy as $\arg \nu$ goes from $\pi/2$ towards $0$ and becomes
completely useless for real $\nu$. A numerical comparison of the
values obtained with Eqs. (\ref{docho}) and (\ref{vuno}) for
moderate values of $s$ and $n$ shows that both labels coincide.

\section{Zeros of $K_{\nu}(z)$ for intermediate values of $|\nu|$.}

We have used a numerical procedure to explore the behaviour of the
$z$-zeros of $K_{\nu}(z)$, as a function of $\nu$, for moderate
values of $|\nu|$. In order to simplify the presentation of the
results, two different situations have been considered: (i)
variable $|\nu|$ with fixed $\arg \nu$, and (ii) variable $\arg
\nu$ keeping $|\nu|$ fixed. The resulting trajectories of the
first, second, and third zeros are represented in Figures 1, 2,
and 3, respectively. For the evaluation of $K_{\nu}(z)$,
expression (\ref{tres}) has been used. More sophisticated
algorithms can be found in the literature
\cite{camp,fab2,gaut,gil1,gil2,gil3,gil4,segu,tem1,tem2,thom}, but
either they have been conceived only for particular (real or pure
imaginary) values of the variable or the order, or they become
unnecessarily complicated for our purpose of giving a merely
qualitative description of the zeros. For the location of these,
the Newton method has been used. Let us mention, however, other
existing techniques like, for instance, the auxiliary tables for
the evaluation of the $\nu$-zeros of $K_{\nu}(z)$, for $z$ real or
pure imaginary, proposed by Cochran and Hoffspiegel \cite{cochof},
the finite element approximation method, applied by Leung and
Ghaderpanah  \cite{leun} to a very precise determination of the
zeros of $K_n(z)$, and a procedure, applied by Segura \cite{seg}
to the computation of zeros of Bessel and other special functions,
that uses fixed point iterations and does not require the
evaluation of the functions.

\begin{figure}
\epsfig{figure=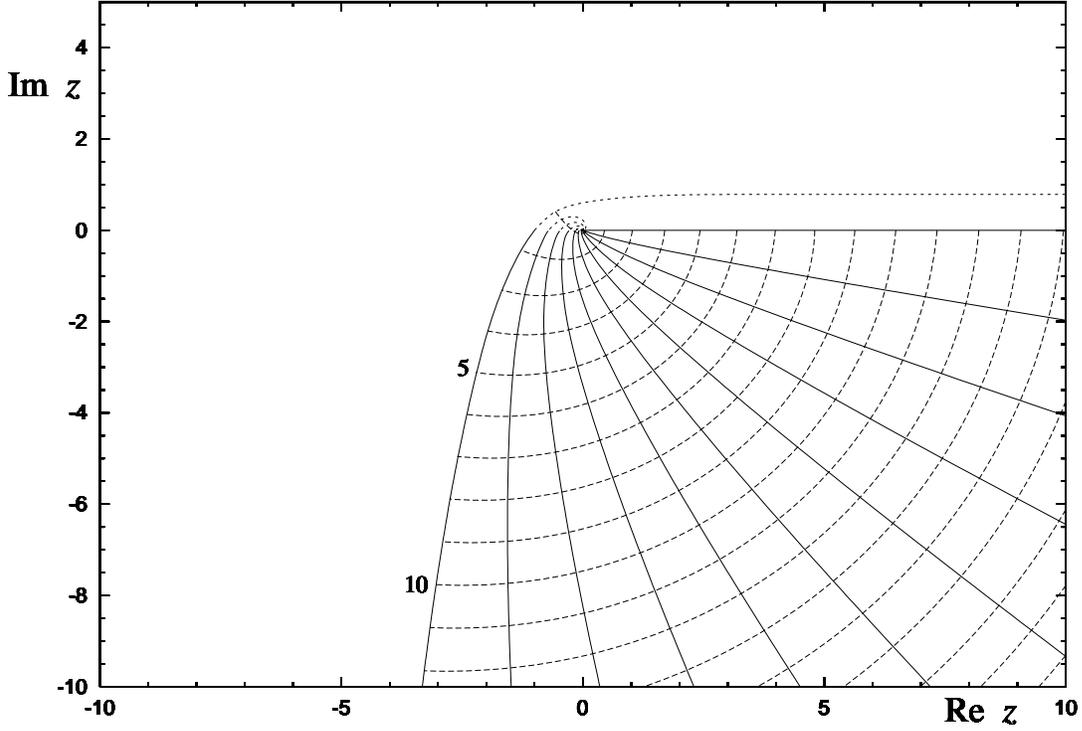,height=10cm,angle=0}
\caption{Trajectories of the first zero, $z_1$, of $K_{\nu}(z)$.
The continuous lines represent the trajectories followed by $z_1$
as $\nu$ varies whereas $\arg \nu$ remains fixed and equal to
$m\pi/20$, with $m=0$ (first trajectory from the left), $1, 2,
\ldots , 10$ (horizontal trajectory). The continuous lines turn
into short-dashed ones as the zero leaves the Riemann sheet
$-\pi<\arg z\leq \pi$. The trajectory corresponding to real $\nu$
($m=0$) goes to $(+\infty, \pi/4)$, with $\arg z_1\to-2\pi$, for
$\nu=1/2$. The nearly circular dashed arcs represent the
trajectories of $z_1$ for fixed integer $|\nu|$ and varying $\arg
\nu$ going from $0$ to $\pi/2$. The trajectories corresponding to
$|\nu|=5$ and $10$ are explicitly indicated.} \label{Macf1}
\end{figure}

\begin{figure}
\epsfig{figure=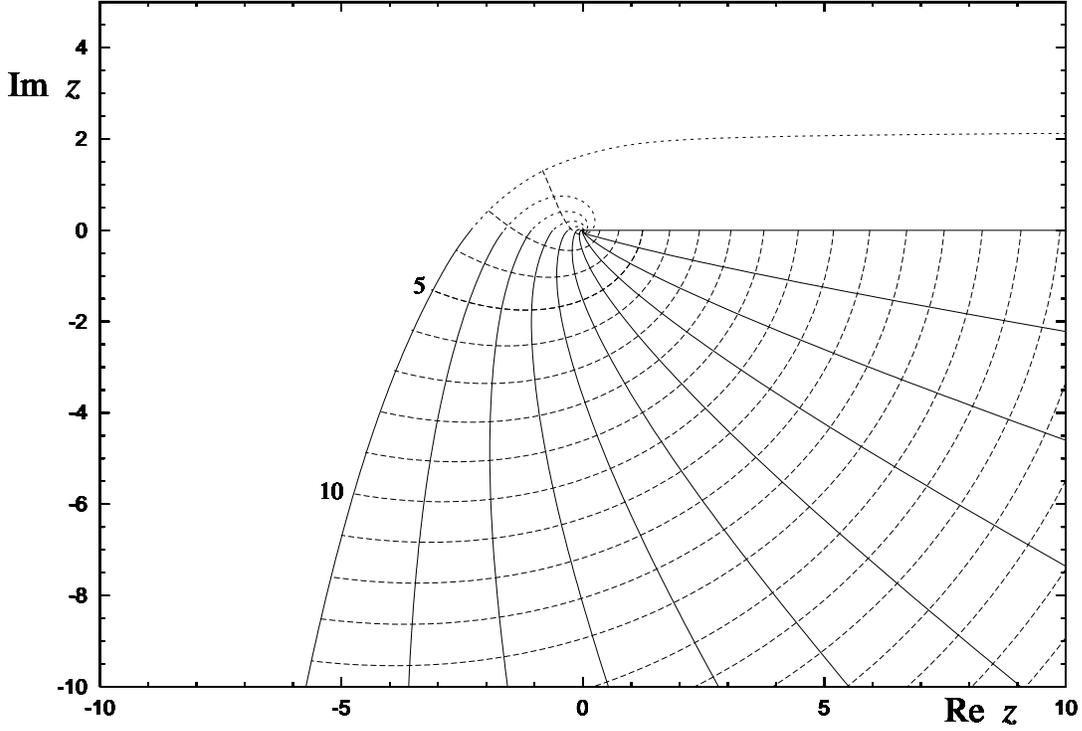,height=10cm,angle=0}
\caption{Trajectories of the second zero, $z_2$, of $K_{\nu}(z)$.
The same convention as in Fig. 1 has been adopted to represent the
behaviour of $z_2$. The trajectory for real $\nu$ goes to
$(+\infty, 3\pi/4)$, with $\arg z_2\to-2\pi$, for $\nu=3/2$.}
\label{Macf2}
\end{figure}

\begin{figure}
\epsfig{figure=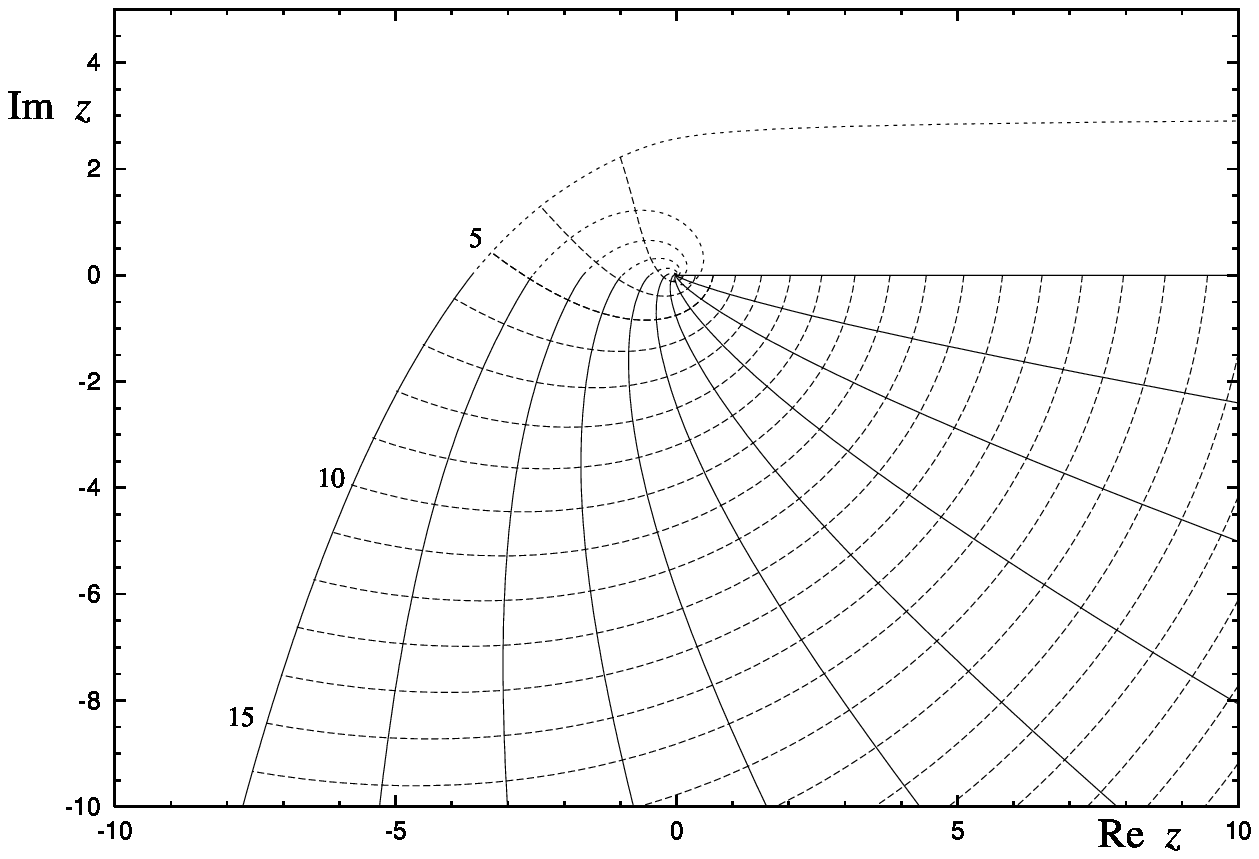,height=10cm,angle=0}
\caption{Trajectories of the third zero, $z_3$, of $K_{\nu}(z)$.
The lines represent the behaviour of $z_3$ with the same
convention as in Fig. 1. The trajectory for real $\nu$ goes to
$(+\infty, 5\pi/4)$, with $\arg z_3\to-2\pi$, for $\nu=5/2$.}
\label{Macf3}
\end{figure}

A double precision FORTRAN code was used to compute the expression
\begin{equation}
\mathcal{E}_{\nu}(z) = \sum_{k=0}^{\infty}
\frac{(z^2/4)^k}{k!(1\!-\!\nu)_k}
-(z/2)^{2\nu}\frac{\Gamma(1\!-\!\nu)}{\Gamma(1\!+\!\nu)}
\sum_{k=0}^{\infty} \frac{(z^2/4)^k}{k!(1\!+\!\nu)_k}
\label{extra}
\end{equation}
whose zeros coincide with those of the right hand side of Eq.
(\ref{tres}). The two series were summed term by term until the
absolute value of the ratio of the last term to the partial sum
turned out to be less than $10^{-18}$. A warning message was
foreseen for the case of the number of summed terms becoming
larger than $100$, but this limit was never reached. Fortunately,
zeros at large values of $|z|$ occur for values of $\nu$ either of
large modulus or in the neighbourhood of $n+\frac{1}{2}$, $n=0, 1,
2, \ldots$. In the first case, the series in (\ref{extra})
converge reasonably, from the numerical point of view; in the
second case, analytical methods can be used to approximate the
zeros. Typical runs of our procedure can be seen in Table 1, which
shows intermediate results in the determination of the first,
second and third zeros of $K_{\nu}(z)$, with $|\nu|=21$, $\arg \nu
= 7\pi/20$. Initial approximate values of the zeros, for a given
$\nu$, were obtained by extrapolation of the trajectories followed
by the zeros as $\nu$ varies. Besides the smallness of the values
of the real and imaginary parts of $\mathcal{E}_{\nu}(z)$, their
changes of sign in the successive steps of the Newton method make
us to be confident about our results.

\begin{table}
\caption{Typical output of the Newton method applied to the
location of the first, second and third zeros of $K_{\nu}(z)$,
with $|\nu|=21$, $\arg \nu = 7\pi/20$. Besides the successive
approximations to the zeros, the values of $\mathcal{E}_{\nu}(z)$,
given by Eq. (\ref{extra}), are shown.}
\begin{tabular}{rrrrrrr}
\hline $\Re z_1\qquad$ & $\;$ & $\Im z_1\qquad$ & $\qquad$ &
$\Re\mathcal{E}_{\nu}(z_1)\quad$ & $\;$ &
$\Im\mathcal{E}_{\nu}(z_1)\quad$ \\
\hline
14.00000000 & & $-$ 8.60000000 &  & $-$ .127330E+01 & & $-$ .817656E$-$01 \\
14.02983390 & & $-$ 8.68222831 &  &   .154771E+00 & &   .778107E$-$01 \\
14.02406924 & & $-$ 8.67520518 &  &   .105223E$-$01 & &   .979004E$-$03 \\
14.02389450 & & $-$ 8.67463545 &  & $-$ .594713E$-$04 & &   .105606E$-$04 \\
14.02389461 & & $-$ 8.67463884 &  &   .429746E$-$07 & & $-$ .440585E$-$08 \\
\hline \hline  $\Re z_2\qquad$ & $\;$ & $\Im z_2\qquad$ & $\qquad$
& $\Re\mathcal{E}_{\nu}(z_2)\quad$ & $\;$ &
$\Im\mathcal{E}_{\nu}(z_2)\quad$ \\
\hline
10.90000000 & & $-$ 8.00000000 &  & $-$ .151232E+01 & & $-$ .835798E+00 \\
10.98486986 & & $-$ 7.94986297 &  & $-$ .232404E+00 & &   .536993E$-$01 \\
11.00073780 & & $-$ 7.95926858 &  &   .207686E$-$01 & & $-$ .296534E$-$01 \\
10.99984823 & & $-$ 7.95696547 &  &   .191900E$-$03 & & $-$ .214639E$-$03 \\
10.99983888 & & $-$ 7.95694790 &  & $-$ .278481E$-$06 & &   .649298E$-$06 \\
10.99983889 & & $-$ 7.95694795 &  &   .698624E$-$11 & & $-$ .156303E$-$10 \\
\hline \hline $\Re z_3\qquad$ & $\;$ & $\Im z_3\qquad$ & $\qquad$
& $\Re\mathcal{E}_{\nu}(z_3)\quad$ & $\;$ &
$\Im\mathcal{E}_{\nu}(z_3)\quad$ \\
\hline
8.50000000 & & $-$ 7.00000000 &  &   .741241E+00 & &   .345969E+01 \\
8.60000000 & & $-$ 6.92107409 &  &   .978987E+00 & &   .293843E+01 \\
8.70000000 & & $-$ 7.02107409 &  &   .658986E+00 & &   .247911E+01 \\
8.80000000 & & $-$ 7.12107409 &  &   .598150E+00 & &   .174818E+01 \\
8.88169635 & & $-$ 7.22107409 &  &   .831925E+00 & &   .871598E+00 \\
8.80046094 & & $-$ 7.31734050 &  & $-$ .282066E+00 & &   .243279E+00 \\
8.82858965 & & $-$ 7.32026067 &  &   .229670E$-$01 & &   .745812E$-$01 \\
8.82891727 & & $-$ 7.32660986 &  &   .257068E$-$04 & & $-$ .695341E$-$03 \\
8.82889674 & & $-$ 7.32655881 &  & $-$ .673649E$-$06 & & $-$ .726223E$-$05 \\
8.82889659 & & $-$ 7.32655825 &  &   .622240E$-$10 & &   .575485E$-$09 \\
\hline \end{tabular} \end{table}

Our figures show that for large $|\nu|$, with $\arg \nu\in (0,
\pi/2)$, the zeros lie far from the origin in the half plane
$-\pi<\arg z<0$, according to Eq. (\ref{docho}). When $|\nu|$
decreases, they approach the origin but, before reaching it, they
cross the semiaxis $\arg z=-\pi$, leaving in this way the
principal Riemann sheet $-\pi<\arg z\leq\pi$, for a certain
$|\nu_{{\rm crit},s}|$ whose value depends on $\arg \nu$ and on
the order $s$ of the zero $z_s$. The smaller $\arg \nu$ (between
$0$ and $\pi/2$) and the higher the label $s$ of the zero, the
larger the value of $|\nu_{{\rm crit},s}|$. (In the case of real
$\nu$, it can be shown that $\nu_{{\rm crit},s}=2s-1/2$.)
Therefore, except for the case of pure imaginary $\nu$, the number
of zeros of $K_{\nu}(z)$ in the principal Riemann sheet is finite:
only those zeros with label $s=1, 2, \ldots$ such that $|\nu_{{\rm
crit},s}|<|\nu|$, for the corresponding value of $\arg \nu$,
remain in that principal Riemann sheet. If $|\nu|$ decreases
again, being $\arg \nu \neq 0$, the zeros approach the origin
spiraling clockwise infinitely, according to Eqs. (\ref{duno}) and
(\ref{ddos}). The case of real $\nu$ deserves a particular
consideration.

\section{Zeros of $K_{\nu}(z)$ for real $\nu$.}

Let us now concentrate on the particular case of $\nu$ being real.
According to the second of Eqs. (\ref{uno}), the zeros of
$K_{\nu}(z)$ appear in complex conjugate pairs. We refer, in what
follows, only to those zeros of negative argument. Let us also
recall that, for half-integer $\nu$, $\nu=n+\frac{1}{2}$, the
Macdonald function reduces to a factor, different from zero in the
finite plane, times a polynomial of degree $n$. Therefore,
$K_{\nu}(z)$ presents exactly $n$ $z$-zeros in each one of the
Riemann sheets. Of course, the locations of these zeros are the
same in all sheets.

The aproximate expression (\ref{docho}) for the first zeros of
$K_{\nu}(z)$ ($\nu\gg 1$) is valid also in this case. Equation
(\ref{vuno}), instead, is not applicable because the assumption
$|z|\ll |\nu|$, used in its derivation, is not justified. Also,
Eqs. (\ref{duno}) and (\ref{ddos}) are no more valid due to the
possibility of $|z|$ becoming infinite for certain values of
$\nu$, a fact that invalidates Eq. (\ref{cinco}), from which Eqs.
(\ref{duno}) and (\ref{ddos}) stem.

According to Eq. (2), the $z$-zeros of $K_{\nu}(z)$ are obtained
from the $z$-zeros of $H_{\nu}^{(1)}(z)$ by a rotation of $-\pi/2$
in the complex $z$-plane. It is then immediate to deduce the
behaviour of the zeros of $K_{\nu}(z)$ for varying real $\nu$ from
a previous discussion \cite{cru1,cru2} of the zeros of
$H_{\nu}^{(1)}(z)$. For convenience of the reader, we recall here,
without demonstration, the main results of that discussion
translated to those zeros in the lower half plane $\Im z<0$. To be
specific, let us focus on the trajectory followed by a zero,
$z_s$, as $\nu$ decreases from $+\infty$ to $0$.

For very large $\nu$, $z_s$ is given approximately by Eq.
(\ref{docho}). It lies far from the origin in the quadrant
$-\pi<\arg z<-\pi/2$. As $\nu$ decreases, it moves upwards and to
the right, reaching the horizontal semiaxis $\arg z = -\pi$ for
$\nu=2s-1/2$. Then, it makes a nearly circular quarter of a turn,
crosses the vertical semiaxis $\arg z = -3\pi/2$ for $\nu=s-1/3$,
and goes nearly horizontally towards $(+\infty, (2s-1)\pi/4)$ as
$\nu\to s-1/2$ from above.

The behaviour of $z_s$ for $\nu<s-1/2$ is not uniquely defined. In
fact, the implicit function $z_s(\nu)$ defined by the condition
$K_{\nu}(z)=0$ presents a logarithmic branch point at $\nu=s-1/2$,
and the values of $z_s$ for $\nu<s-1/2$ depend on the chosen
branch. One can pass continuously from the values of $z_s$ for
$\nu>s-1/2$ to those for $\nu<s-1/2$ if one avoids the branch
point by the usual procedure of adding a small imaginary part to
$\nu$ as its real part crosses the value $s-1/2$. But the
resulting values of $z_s$ for $\nu<s-1/2$ depend on the sign of
that imaginary part. Since this ambiguity in the values of $z_s$
has been thoroughly discussed in Ref. \cite{cru2}, in the context
of the $z$-zeros of $H_{\nu}^{(1)}(z)$, we do not consider
necessary to pursue further.

Besides those described above, $K_{\nu}(z)$ presents an infinite
set of zeros near and  at the left of the vertical semiaxis $\arg
z =-3\pi/2$ for every integer $\nu=n$ \cite[p. 513]{wats}. As
$\nu$ increases, they approach that semiaxis, cross it when
$\nu=n+1/3$, and move nearly horizontally to the right, going to
infinity as $\nu\to n+1/2$, which is a logarithmic branch point
for the position of the zeros as a function of $\nu$, as it was
shown  in Refs. \cite{cru1,cru2} for the zeros of
$H_{\nu}^{(1)}(z)$. The zeros make a discontinuous jump of $\pm
i\pi/2$, according to how the branch point is circumvented, move
nearly horizontally to the left, reach the vertical semiaxis $\arg
z =-3\pi/2$ when $\nu=n+2/3$, and go to the left, to end, for
$\nu=n+1$, at positions intermediate between those that they
occupied for $\nu=n$. A figure illustrating that behaviour can by
obtained by rotating Fig. 1 in ref. \cite{cru1} by an angle of
$-\pi/2$.

\section*{Acknowledgements}

One of the authors (EF) acknowledges and thanks support received
from Conselho Nacional de Desenvolvimento Cient\'{\i}fico e
Tecnol\'ogico (CNPq, Brazil). The other author (JS) is grateful to
Dr. Yik-Man Chiang for illuminating discussions and acknowledges
financial aid of Comisi\'on Interministerial de Ciencia y
Tecnolog\'{\i}a and of Diputaci\'on General de Arag\'on.


\begin{thebibliography}{99}

\bibitem{abra} M. Abramowitz and I. A. Stegun (Eds.),
Handbook of Mathematical Functions, Dover, New York, 1965.

\bibitem{camp} J. B. Campbell, On Temme's Algorithm for the
Modified Bessel Function of the Third Kind, ACM Trans. Math.
Softw. 6 (1980) 581--586.

\bibitem{coch} J. A. Cochran, The Zeros of Hankel Functions as
Functions of Their Order, Numer. Math. 7 (1965) 238--250.

\bibitem{cochof} J. A. Cochran and J. N. Hoffspiegel, Numerical
Techniques for Finding $\nu$-Zeros of Hankel Functions, Math.
Comp. 24 (1970) 413--422.

\bibitem{cru1} A. Cruz and J. Sesma, Zeros of the Hankel Function
of Real Order and of Its Derivative, Math. Comp. 39 (1982)
639--645.

\bibitem{cru2} A. Cruz, J. Esparza and J. Sesma, Zeros of the
Hankel function of real order out of the principal Riemann sheet,
J. Comput. Appl. Math. 37 (1991) 89--99.

\bibitem{duns} T. M. Dunster, Bessel functions of purely imaginary
order, with an application to second-order linear differential
equations having a large parameter, SIAM J. Math. Anal. 21 (1990)
995--1018.

\bibitem{fab1} B. R. Fabijonas and F. W. J. Olver, On the
Reversion of an Asymptotic Expansion and the Zeros of the Airy
Functions, SIAM Rev. 41 (1999) 762--773.

\bibitem{fab2} B. R. Fabijonas, D. W. Lozier, J. M. Rappoport,
Algorithms and codes for the Macdonald function: recent progress
and comparisons, J. Comput. Appl. Math. 161 (2003) 179--192.

\bibitem{ferr} E. M. Ferreira and J. Sesma, Zeros of the Modified
Hankel Function, Numer. Math. 16 (1970) 278--284.

\bibitem{fran} W. Franz and R. Galle, Semiasymptotische Reihen
f\"ur die Beugung einer ebenen Welle am Zylinder, Z. Naturforschg.
10a (1955) 374--378.

\bibitem{frie} F. G. Friedlander, Diffraction of Pulses by a Circular
Cylinder, Comm. Pure Appl. Math. 7 (1954) 705--732.

\bibitem{gaut} W. Gautschi, Numerical quadrature computation of
the Macdonald function for complex orders, BIT Numerical
Mathematics 45 (2005) 593--603.

\bibitem{gil1} A. Gil, J. Segura, and N. M. Temme, Evaluation of
the Modified Bessel Function of the Third Kind of Imaginary
Orders, J. Comput. Phys. 175 (2002) 398--411.

\bibitem{gil2} A. Gil, J. Segura, and N. M. Temme, Computing
special functions by using quadrature rules, Numer. Algorithms 33
(2003) 265--275.

\bibitem{gil3} A. Gil, J. Segura, and N. M. Temme, Computing
solutions of the the modified Bessel differential equation for
imaginary orders and positive arguments, ACM Trans. Math. Softw.
30 (2004) 145--158.

\bibitem{gil4} A. Gil, J. Segura, and N. M. Temme, Algorithm 831:
Modified Bessel functions of imaginary order and positive
argument, ACM Trans. Math. Softw. 30 (2004) 159--164.


\bibitem{gray} A. Gray and G. B. Mathews, A Treatise on Bessel
Functions and Their Application to Physics, 2nd edition, Dover,
New York, 1966.

\bibitem{heth} H. W. Hethcote, Error Bounds for Asymptotic
Approximations of Zeros of Hankel Functions Occurring in
Diffraction Problems, J. Math. Phys. 11 (1970) 2501--2504.

\bibitem{kell} J. B. Keller, S. I. Rubinow, and M. Goldstein,
Zeros of Hankel Functions and Poles of Scattering Amplitudes, J.
Math. Phys. 4 (1963) 829--832.

\bibitem{lafo} A. Laforgia, Inequalities and monotonicity results
for zeros of modified Bessel functions of purely imaginary order,
Quart. Appl. Math. 44 (1986) 91--96.

\bibitem{leun} K. V. Leung and S. S. Ghaderpanah, An Application
of the Finite Element Approximation Method to find the Complex
Zeros of the Modified Bessel Function $K_n(z)$, Math. Comp. 33
(1979) 1299--1306. Notice that a factor $10^{-1}$ is lacking in
the values given in this reference for the imaginary parts of the
first zero of $K_2(z)$, the second one of $K_4(z)$, the third one
of $K_6(z)$, the fourth one of $K_8(z)$, and the fifth one of
$K_{10}(z)$.

\bibitem{loz1} D. W. Lozier, Software needs in special functions,
J. Comput. Appl. Math. 66 (1996) 345--358.

\bibitem{loz2} D. W. Lozier and F. W. J. Olver, Numerical
evaluation of special functions, in: W. Gautschi, Ed., Mathematics
of Computation 1943--1993: A Half-Century of Computational
Mathematics, Proc. Symposia in Applied Mathematics, Vol 48,
American Mathematical Society, Providence, RI, 1994, pp. 79--125.

\bibitem{magn} W. Magnus and L. Kotin, The Zeros of the Hankel
Function as a Function of its Order, Numer. Math. 2 (1960)
228--244.

\bibitem{nale} G. F. Nalesso, On the Zeros of a Class of Bessel
Functions Whose Argument and Order are Functions of a Complex
Variable, IMA J. Appl. Math. 43 (1989) 195--217.

\bibitem{sand} S. E. Sandstr\"om and C. Ackr\'en, Note on the
complex zeros of $H_{\nu}^{\prime}(x)+i\zeta H_{\nu}(x)\!=\!0$. J.
Comput. Appl. Math., in press. Available online,
DOI:10.1016/j.cam.2006.01.032.

\bibitem{seg} J. Segura, The zeros of special functions from a
fixed point method, SIAM J. Numer. Anal. 40 (2002) 114-133.

\bibitem{segu} J. Segura, P. Fern\'andez de C\'ordoba, Yu. L.
Ratis, A code to evaluate modified Bessel functions based on the
continued fraction method, Comput. Phys. Commun. 105 (1997)
263--272.

\bibitem{tem1} N. M. Temme, On the Numerical Evaluation of the
Modified Besel Function of the Third Kind, J. Comput. Phys. 19
(1975) 324--337.

\bibitem{tem2} N. M. Temme, Numerical algorithms for uniform
Airy-type asymptotic expansions, Numer. Algorithms 15 (1997)
207--225.

\bibitem{thom} I. J. Thompson and A. R. Barnett, Modified Bessel
functions $I_{\nu}(z)$ and $K_{\nu}(z)$ of real order and complex
argument, to selected accuracy, Comput. Phys. Commun. 47 (1987)
245--257.

\bibitem{wats} G. N. Watson, A Treatise on the Theory of Bessel
Functions, 2nd edition, Cambridge Univ. Press, Cambridge, 1944.

\end{thebibliography}
\end{document}